\documentclass[12pt,francais]{smfart}
\usepackage[T1]{fontenc}
\usepackage[english,francais]{babel}
\usepackage{amssymb,mathrsfs}
\usepackage[matrix,arrow,curve,cmtip]{xy}

 \calclayout

\def\Spec{\operatorname{Spec}}
\def\Div{\operatorname{Div}}

\def\Sym{\operatorname{Sym}}
\def\div{\operatorname{div}}
\def\fp{\mathfrak{p}}
\def\gm{{\mathbf{G}_m}}
\def\hPic{\widehat{\operatorname{Pic}}\vphantom{\operatorname{Pic}}}
\def\hdeg{\widehat{\operatorname{deg}}\,}
\def\hdiv{\widehat{\operatorname{div}}\,}
\def\id{\mathrm{id}}
\def\mcL{\mathscr L}
\def\sbullet{{\scriptscriptstyle\bullet}}
\def\Proj{\operatorname{Proj}}
\def\Pic{\operatorname{Pic}}
\def\V{{\mathbf V}}
\def\Z{{\mathbf Z}}
\def\P{{\mathbf P}}
\def\C{{\mathbf C}}
\def\Q{{\mathbf Q}}
\def\R{{\mathbf R}}
\def\norm#1{\left\|{#1}\right\|}
\def\Hom{\operatorname{Hom}}
\let\ra\rightarrow \let\hra\hookrightarrow
\let\union\cup
\let\eps\varepsilon \let\phi\varphi

\def\qedsymbol{\ensuremath{\blacksquare}}

\theoremstyle{plain} 
\newtheorem{defi}{D\'efinition}[section]
\newtheorem{prop}[defi]{Proposition}
\newtheorem{theo}[defi]{Th\'eor\`eme}
\newtheorem{coro}[defi]{Corollaire}
\newtheorem{lemm}[defi]{Lemme}

\theoremstyle{remark}
\newtheorem{rema}[defi]{Remarque}

\begin{document}

\title[Hauteurs canoniques sur des vari\'et\'es semi-ab\'eliennes]
      {G\'eom\'etrie d'Arakelov \\ et hauteurs canoniques \\
       sur des vari\'et\'es semi-ab\'eliennes}
\alttitle{Canonical heights and Arakelov geometry on semi-abelian varieties}
\author{Antoine Chambert-Loir}
\address{Institut de math\'ematiques de Jussieu\\ Boite 247 \\
4, place Jussieu \\ F-75252 Paris Cedex 05}
\email{chambert@math.jussieu.fr}


\keywords {Arakelov geometry, theorem of the square, semi-abelian
variety, relative height.}
\subjclass{14G, 14K, 11G, 14G40, 14K15}

\begin{abstract}
Cet article propose une construction en th\'eorie d'Arakelov
des hauteurs canoniques sur une extension d'une vari\'et\'e ab\'elienne
par le groupe multiplicatif. Celles-ci apparaissent comme la somme
d'une hauteur provenant de la vari\'et\'e ab\'elienne et
de ce que nous appelons {\em hauteur relative.} 
Nous \'etudions aussi les points de hauteur relative nulle.
\end{abstract}

\begin{altabstract}
\selectlanguage{english}
In this paper, we propose a construction of the canonical heights
on an extension of an abelian variety by the multiplicative group,
in the framework of Arakelov geometry. These canonical heights
are the sum of some height coming from the abelian variety and
something we call a {\em relative height.}
We finally give some complements about the points whose relative
height is zero.
\end{altabstract}

\maketitle

\setcounter{tocdepth}{2}
\tableofcontents

\section{Introduction}

Dans cet article, nous voulons montrer comment la th\'eorie d'Arakelov
permet d'interpr\'eter les  hauteurs canoniques au sens
de~\cite{bertrand96,call-s93}
sur une extension d'une vari\'et\'e ab\'elienne par 
le groupe multiplicatif $\gm$. Dans l'esprit de la construction
arakelovienne de la hauteur de N\'eron--Tate 
(cf.~\cite{moret-bailly85,moret-bailly85b,faltings-w85}), 
nous montrons que sur une
telle extension, il existe une hauteur canonique et elle est donn\'ee par
le degr\'e d'Arakelov d'un fibr\'e inversible 
sur un mod\`ele convenable muni de m\'etriques hermitiennes
aux places archim\'ediennes. Le cas d'une
vari\'et\'e semi-ab\'elienne dont le tore sous-jacent est {\em
d\'eploy\'e\/} se traite par les m\^emes m\'ethodes, voir la
remarque~\ref{rema:semiab}.

Le mod\`ele entier est donn\'e gr\^ace \`a 
la formule de Weil--Barsotti dans le cas de bonne r\'eduction
et \`a une extension de cette formule faisant intervenir
la composante neutre des mod\`eles de N\'eron en g\'en\'eral
(cf.~\cite[(5.1), p.~53]{mazur-m76}). Nous rappelons ceci au
paragraphe~\ref{sec:t-barsotti}.

Sur une variante enti\`ere de la compactification de Serre,
Faltings--W\"ustholz~\cite{serre79,faltings-w86}
que nous exhibons au paragraphe~\ref{sec:t-comp},
nous produisons des faisceaux inversibles relativement amples
et les munissons de m\'etriques hermitiennes \`a l'infini.
Ceci fait,
nous montrons au paragraphe~\ref{sec:t-haut}
que l'on obtient la hauteur canonique en calculant le degr\'e d'Arakelov
d'un des fibr\'es inversibles m\'etris\'es pr\'ec\'edemment d\'efinis
(th\'eor\`eme~\ref{theo:t-hcan} et corollaire~\ref{coro:t-hcan}).
La preuve est alors analogue \`a celle de~\cite{moret-bailly85b,moret-bailly85} :
le manque d'uniformit\'e des mod\`eles entiers est compens\'e par les
propri\'et\'es du degr\'e d'Arakelov calcul\'e relativement aux morphismes
de multiplication par un entier sur le groupe alg\'ebrique.

Nous donnons au paragraphe~\ref{sec:t-misc}
contient quelques compl\'ements sur les
points de hauteur relative nulle et les {\og points de Ribet\fg}
de~\cite{jacquinot-r87,bertrand96}.
Dans le cas de bonne r\'eduction, nous obtenons la caract\'erisation suivante 
(proposition~\ref{prop:rel-0}) :

{\itshape 
Soient $K$ un corps de nombres et $\mathfrak o_K$ l'anneau des
entiers de $K$. Soient $\mathscr A$ un $\mathfrak o_K$-sch\'ema ab\'elien
et $1\ra \gm \ra \mathscr E\ra\mathscr A \ra 0$ une extension de
$\mathscr A$ par $\gm$ fournie par un faisceau inversible $\mathscr
L\in\Pic^0(\mathscr A)$. La section nulle de $\mathscr E$
induit une rigidification de $\mathscr L$ \`a l'origine de $\mathscr A$,
laquelle rigidification d\'etermine un isomorphisme du carr\'e.

Le faisceau inversible $\mathscr L\otimes_{\mathfrak o_K}\C $
sur $\mathscr A\otimes_{\mathfrak o_K}\C$
admet alors une unique m\'etrique hermitienne 
telle que l'isomorphisme du carr\'e soit une isom\'etrie. 

Soit $x\in\mathscr A_K(K)$ et $\eps_x:\Spec\mathfrak o_K\ra\mathscr A$
l'unique section qui prolonge $x$. Il existe alors un point de hauteur
relative nulle dans $\mathscr E_K(K)$ au-dessus de $x$ si et
seulement si l'\'el\'ement $\eps_x^*\mathscr L$ de $\hPic(\Spec\mathfrak o_K)$
est trivial, c'est-\`a-dire admet une base de norme~$1$ en toute place.
}

Nous terminons cet article en \'evoquant bri\`evement
comment l'on peut le formuler dans le langage des m\'etriques
ad\'eliques de S.~Zhang.

\bigskip

{\em\small
Cet article est une version l\'eg\`erement remani\'ee du
premier chapitre de ma th\`ese~\cite{chambert95}, soutenue en
d\'ecembre 1995\dots

Je tiens \`a remercier
chaleureusement mon directeur de th\`ese, Daniel Bertrand,
pour son aide et ses encouragements incessants durant
la gestation de ce travail.
Je remercie enfin Jean-Beno\^\i t Bost et Ahmed Abbes pour leurs
remarques.

}

\section{Notations et conventions} 

Si $X$ est un espace localement annel\'e
et $\mathscr F$ un faisceau quasi-coh\'erent sur~$X$, on utilise les conventions
de [EGA II] en notant $\V(\mathscr F)=\Spec\Sym^\sbullet \mathscr F$ et
$\P(\mathscr F)=\Proj\Sym^\sbullet\mathscr F$ les fibr\'es {\og
vectoriels\fg}
et {\og projectifs\fg} associ\'es \`a $\mathscr F$. En particulier, un morphisme
$u:\mathscr E\ra\mathscr F$ d\'efinit des applications dans l'autre sens
$\V(\mathscr F)\ra\V(\mathscr E)$
et $\P(\mathscr F)\ra\P(\mathscr E)$, cette derni\`ere n'\'etant 
d\'efinie que sur un ouvert si $u$ n'est pas surjective.
Nous commettrons l'abus de langage consistant \`a appeler fibr\'e en droites
un faisceau localement libre de rang~1.

\medskip

Soit $X$ un sch\'ema plat et quasi-projectif sur $\Z$. Un fibr\'e
en droites m\'etris\'e sur $X$ est la donn\'ee d'un fibr\'e en droites
$\mathscr L$ sur $X$, ainsi que d'une m\'etrique hermitienne (continue)
sur le fibr\'e complexe $\mathscr L_\C$ sur $X(\C)$.  On 
demandera que la m\'etrique hermitienne soit compatible \`a la conjugaison 
complexe.
On note alors $\hPic(X)$ le groupe ab\'elien pour le produit tensoriel
des classes d'isomorphisme de fibr\'es en droites m\'etris\'es.
Tout morphisme de sch\'emas $f:X\ra X'$ induit un morphisme
de groupes $f^*:\hPic(X')\ra\hPic(X)$.

Soit $K$ un corps de nombres, $\mathfrak o_K$ son anneau d'entiers
et notons $S=\Spec\mathfrak o_K$.
Les \'el\'ements de $\hPic(S)$ sont alors les (classes
d'isomorphisme de)  $\mathfrak o_K$-modules
projectifs $\mathscr L$ de rang~$1$ munis d'une m\'etrique hermitienne 
sur les droites complexes $\mathscr 
L\otimes_\sigma\C$ (compatibles \`a la conjugaison 
complexe).
Un \'el\'ement de $\hPic(S)$ poss\`ede un degr\'e d'Arakelov, d\'efini
par la formule
$$ \hdeg (\mathscr L,\|\cdot\|_\sigma) = \log \# (\mathscr L/s\mathfrak 
o_K)  - \sum_{\sigma:K\hra\C} \log \|s\|_\sigma, $$
o\`u $s$ est un \'el\'ement non nul quelconque de $\mathscr L$~; d'apr\`es 
la formule du produit, il est ind\'ependant du choix de $s$.
L'application $\hdeg:\hPic(S)\ra\R$ est un homomorphisme
de groupes ab\'eliens.

\medskip

Soit  $X$ un sch\'ema plat et projectif sur l'anneau des 
entiers d'un corps de nombres $K$ et $(\mathscr L,\|\cdot\|)\in \hPic(X)$.
Associons \`a tout point $P\in X_K(\overline K)$ le r\'eel
$h(P)=[K':\Q]^{-1}\hdeg \eps_P^*(\mathscr L,\|\cdot\|)$ o\`u $K'$ est un corps
de d\'efinition de $P$ et $\eps_P:\Spec \mathfrak o_{K'}\ra X$ est la 
section canonique. Alors, la fonction $P\mapsto h(P)$ est un repr\'esentant
de la hauteur de Weil (logarithmique, absolue)
de~$P$ pour le fibr\'e en droites $\mathscr L_K$
sur $X_K$.
(Voir~\cite{szpiro85}, ou~\cite{bost-g-s94} pour des g\'en\'eralisations.)

\section{Formule de Weil--Barsotti}
\label{sec:t-barsotti}
Commen\c{c}ons par rappeler cette formule dans le cas de sch\'emas
ab\'eliens. Soient $S$ un sch\'ema noeth\'erien et $A$ un $S$-sch\'ema
ab\'elien. D'apr\`es Hilbert~90, il correspond \`a
une $S$-extension commutative $1\ra\gm\ra E\ra A\ra 0$
de $A$ par $\gm$ un espace principal homog\`ene
sous $\gm$ sur $A$ et donc un faisceau inversible $\mcL\in \Pic(A)$
tel que $E$ s'identifie au fibr\'e en droites $\V(\mcL^\vee)$ 
priv\'e de sa section nulle. La section neutre $\eps_E:S\ra E$ correspond
\`a un isomorphisme $\mathscr O_S\simeq \eps_A^*\mathscr L$,
c'est-\`a-dire \`a une rigidification de $\mathscr L$ le long de
la section neutre de $A$.
Notons $m$ (resp.\ $p_1$, $p_2$) l'addition (resp. les deux projections)
$A\times_S A \ra A$. 

\begin{prop}[Barsotti--Rosenlicht--Weil, \protect\cite{raynaud68}] \label{prop:wb}
Soient $S$ un sch\'ema et 
$A$ un sch\'ema ab\'elien sur $S$. On note $A^\vee=\Pic^0(A/S)$
le sch\'ema ab\'elien dual.
L'application $E\mapsto\mcL$ d\'ecrite ci-dessus est un isomorphisme 
de foncteurs en groupes sur la cat\'egorie des $S$-sch\'emas
$$ \mathbf{Ext}^1_S(A,\gm) \xrightarrow{\sim} A^\vee \ . $$
\end{prop}
\begin{proof} (cf.~\cite[Appendice]{moret-bailly81}).
Montrons d'abord que cette application est bien \`a valeurs dans
$A^\vee$. En effet, si $S'$ est un $S$-sch\'ema, $x\in A(S')$ et $\xi\in E(S')$
rel\`eve $x$, la translation $T_x$ par $x$ dans $A_{S'}$ 
(resp. $T_\xi$ par $\xi$ dans $E_{S'}$) nous fournit un diagramme commutatif
$$  
\xymatrix{
 {E_{S'}} \ar^{T_\xi}[r] \ar[d] & {E_{S'}} \ar[d] \\
 {A_{S'}} \ar^{T_x}[r] & {A_{S'}\rlap{\quad ,}} }
$$
d'o\`u un isomorphisme $\mcL \ra T_x^* \mcL$.


Cela implique alors que le fibr\'e inversible 
$m^*\mcL \otimes p_1^*\mcL^{-1} \otimes p_2^*\mcL^{-1}$
sur $A\times_S A$ est trivial et donc que $\mcL\in\Pic^0(A)(S)=A^\vee(S)$
(cf.~\cite{raynaud68}).

Il est imm\'ediat que cette application est un morphisme de groupes.
Montrons qu'elle est injective, autrement dit qu'il
existe une unique structure d'extension de $A$ par $\gm$ sur
le sch\'ema~$E_0=\gm\times_S A$. En effet, la multiplication dans $E_0$
s'interpr\`ete comme une application  $A\times_S A \ra \gm$ qui
est n\'ecessairement constante ($A\times_S A$ est propre sur $S$, tandis
que $\gm_S$ est affine) donc nulle, si bien que l'extension
consid\'er\'ee est triviale.

Enfin, construisons la r\'eciproque de cette application. Soit ainsi
$\mcL\in \Pic^0(A)(S)$ que l'on voit comme un \'el\'ement de $\Pic(A)$
muni d'une rigidification le long de la section neutre,
et posons
$E=\V(\mcL^\vee)\setminus \{0\}$ le fibr\'e
en droites associ\'e priv\'e de sa section nulle.
Le choix de l'\'el\'ement neutre dans $E$ revient \`a se donner
un $S$-point de $E$ au-dessus de $\eps_A$, la section unit\'e de $A$~;
ainsi choisissons une {\og rigidification\fg}
de $\eps_A^* \mcL$ c'est-\`a-dire un isomorphisme 
$\eps_A^*\mcL \xrightarrow{\sim} \mathscr O_S$.
D'autre part, comme $\mcL\in\Pic^0(A)$, il existe un unique
isomorphisme
$m^*\mcL\simeq p_1^*\mcL\otimes p_2^*\mcL$ qui
est compatible avec la rigidification de $\mcL$, d'o\`u
une application
$$ m_E:\V(\mcL^\vee) \times_S 
	\V(\mcL^\vee) \ra \V(\mcL^\vee) $$
qui rel\`eve la multiplication $m:A\times_S A \ra A$ et
compatible avec la section unit\'e $\eps_E:S\ra E$ de~$E$.
C'est la loi de groupe sur $E$ que l'on cherchait. En effet,
l'associativit\'e r\'esulte du fait que les deux compositions
$$
\V(\mcL^\vee)\times_S \left( \V(\mcL^\vee)\times_S \V(\mcL^\vee) \right) 
	\ra \V(\mcL^\vee) $$
et 
$$
\left( \V(\mcL^\vee)\times_S\V(\mcL^\vee) \right) \times_S \V(\mcL^\vee) 
	\ra \V(\mcL^\vee) $$
proviennent toutes deux de l'unique isomorphisme rigidifi\'e
$$ p_{123}^*\mcL \ra p_1^*\mcL \otimes p_2^*\mcL \otimes p_3^*\mcL, $$
$p_{1}$, $p_2$, $p_3$ d\'esignant les projections $A^3\ra A$
et $p_{123}:A^3\ra A$ \'etant l'addition des trois composantes.

De m\^eme, la commutativit\'e de la loi de groupe est une cons\'equence de la 
sym\'etrie de l'isomorphisme rigidifi\'e $m^*\mcL\simeq p_1^*\mcL\otimes 
p_2^*\mcL$.

Enfin, il existe une unique application 
$$ \iota_E :\V(\mcL^\vee)\setminus \{0\} \ra \V(\mcL^\vee)\setminus \{0\} $$
au-dessus de la multiplication par $-1$, compatible avec $\eps_E$
et provenant de la composition de l'isomorphisme 
$[-1]_A^*\mcL\simeq \mcL^\vee$ et de l'application
naturelle 
$\V(\mcL^\vee)\setminus \{0\} \ra \V(\mcL)\setminus\{0\}$
qui associe \`a une base de $\mcL$ la base duale de $\mcL^\vee$.
La compos\'ee $m_E \circ (\id_E,\iota_E)$ est une application $E\ra E$
constante sur les fibres de la projection $E\ra A$ et \`a valeurs dans la
fibre de $E$ au-dessus de $\eps_A$. Elle est ainsi constante et vaut
$\eps_E$, ce qui prouve que $\iota_E$ est le morphisme {\og inverse\fg}.

Nous avons ainsi associ\'e \`a tout \'el\'ement de $\mcL\in\Pic^0(A)(S)$
une extension de $A$ par $\gm$
dont l'espace sous-jacent est $\V(\mcL^\vee)\setminus\{0\}$~;
cette application est la r\'eciproque voulue.
\end{proof}

Soient $S$ un sch\'ema de Dedekind, c'est-\`a-dire
un sch\'ema normal noeth\'erien de dimension~1 et $\pi:A\ra S$ un 
{\og mod\`ele de N\'eron\fg}.
Autrement dit, il existe un ouvert dense $U$ de $S$ tel que $A_U$
est un sch\'ema ab\'elien et $A$ est le mod\`ele de N\'eron de $A_U$ sur $S$.
On note $A^0$ la composante neutre de $A$ c'est-\`a-dire le plus grand
sous-sch\'ema en groupes ouvert de $A/S$ \`a fibres connexes.

Soit $A^\vee$ le mod\`ele de N\'eron dual, c'est-\`a-dire le mod\`ele de
N\'eron du sch\'ema ab\'elien dual $(A_U)^\vee$ (ind\'ependant de $U$). 

On a alors le lemme~:
\begin{prop}[Artin--Mazur, {\cite[Lemme (5.1), p.~53]{mazur-m76}}]
\label{prop:am}
Avec ces notations,
il existe un unique isomorphisme de foncteurs sur la cat\'egorie
des  $S$-sch\'emas lisses
$$ \mathbf{Ext}^1_S(A^0,\gm) \xrightarrow{\sim} A^\vee $$
qui prolonge la dualit\'e des sch\'emas ab\'eliens $A_U$ et $(A_U)^\vee$.
\end{prop}
\begin{proof}
Artin--Mazur prouvent ce lemme en montrant que $\mathbf{Ext}^1_S(A^0,\gm)$
v\'erifie la propri\'et\'e universelle du mod\`ele de N\'eron, \`a savoir
que pour tout $S$-sch\'ema lisse $S'$, une $S'$-extension de $A_{U\times_S
S'}$ par $\gm_{S'}$
se prolonge uniquement en une $S'$-extension de $A^0_{S'}$ par $\gm_{S'}$.
Nous donnons ici une d\'emonstration constructive de la bijectivit\'e
de l'application induite au niveau des $S$-points.

La d\'emonstration de la proposition~\ref{prop:wb} nous ram\`ene \`a
montrer le fait suivant~: soit $\mcL_U \in\Pic(A_U)$ muni d'une
rigidification $\eps_{A_U}^*\mcL_U \simeq \mathscr O_U$ et d'un
isomorphisme rigidifi\'e $m_{A_U}^*\mcL_U \ra p_1^*\mcL_U \otimes p_2^*\mcL_U$,
alors il existe une unique fa\c{c}on de prolonger ces
donn\'ees en un faisceau inversible $\mcL\in\Pic(A^0)$, une rigidification
$\eps_{A^0}^*\mcL \simeq \mathscr O_S$ et un isomorphisme rigidifi\'e
$m_{A^0}^*\mcL \ra p_1^*\mcL \otimes p_2^*\mcL $. (Rappelons
qu'un mod\`ele de N\'eron v\'erifie toujours $\pi_*\mathscr
O_{A^0}=\mathscr O_S$.)

Choisissons un diviseur $D_U\in\Div(A_U)$ tel que $\mcL=\mathscr O(D_U)$.
Comme $A^0$ est r\'egulier, l'adh\'erence sch\'ematique $D$ de $D_U$
dans $A^0$ est un diviseur de $A^0$ et d\'efinissons $\mcL_1=\mathscr O(D)$.
Pour rigidifier $\mcL_1$, on le remplace par
$\mcL_1\otimes \pi^*\eps_{A^0}^*\mcL_1^\vee$ qui est un faisceau inversible
$\mcL_2$ sur $A^0$, muni d'une rigidification qui prolonge la
rigidification initiale sur $\mcL_U$.
Enfin, le faisceau inversible 
$m_{A^0}^*\mcL_2\otimes p_1^*\mcL_2^\vee \otimes p_2^*\mcL_2^\vee$
est (sur $S$) g\'en\'eriquement trivial, puisque trivial une fois restreint \`a
$A_U\times_U A_U \subset A^0\times_S A^0$.
Comme la projection $\pi_2:A^0 \times_S A^0 \ra S$ est \`a fibres connexes
et comme $S$ est de Dedekind, il provient de la base
et est donc de la forme $\pi_2^*\mathscr M$, le faisceau
inversible $\mathscr M_U$ \'etant canoniquement trivial.
Posons finalement $\mcL= \mcL_2 \otimes \pi^*\mathscr M^\vee$. C'est
un \'el\'ement de $\Pic(A^0)$ muni d'une rigidification et d'un
isomorphisme rigidifi\'e comme on voulait, ce qui prouve l'existence
du prolongement.

L'unicit\'e du prolongement se d\'emontre de m\^eme, si $\mcL$ et $\mcL'$
sont deux prolongements, le faisceau inversible $\mcL'\otimes \mcL^\vee$
est g\'en\'eriquement trivial et rigidifi\'e. Il provient ainsi de la base,
mais le faisceau inversible sur $S$ dont il provient est n\'ecessairement
trivial \`a cause des rigidifications. Ainsi, il existe un unique
isomorphisme rigidifi\'e $\mcL\ra \mcL'$ et il est compatible aux deux
isomorphismes suppl\'ementaires.
\end{proof}

\section{Compactification, m\'etriques}
\label{sec:t-comp}
Soit $S$ un sch\'ema de Dedekind connexe, notons
$\eta$ son point g\'en\'erique.
Soient $A_\eta$ une $\eta$-vari\'et\'e ab\'elienne, 
$A$ son mod\`ele de N\'eron sur
$S$ et $A^0$ la composante neutre de $A$.
Soient $E_\eta$ une extension de $A_\eta$ par $\gm$ et $E$ l'extension de $A^0$
par $\gm$ fournie par la proposition~\ref{prop:am}.
Notons  $\mcL$ le faisceau inversible sur $A^0$ associ\'e \`a $E$, vue
comme $\gm$-torseur, de sorte que $E$ s'identifie \`a
$\V(\mcL^\vee)\setminus\{0\}$.

On pose $\mathscr W=\mathscr O_{A^0}\oplus\mcL^\vee$ et
on d\'efinit $\overline E$ comme $\P(\mathscr W)$.
C'est un $A^0$-fibr\'e projectif dont $E$ est un ouvert.
En effet, si $P\in A^0(S)$ et $\eps_P:S\ra A^0$ est la section 
correspondante, un $S$-point de $\overline E$ au-dessus de $P$ correspond \`a
un quotient localement libre de rang~$1$~: $(\alpha,\beta):(\mathscr
O_S\oplus\eps_P^*\mcL^\vee)\ra\mathscr J$. Parmi ceux-ci, les points
de $E$ correspondent aux couples $(\alpha,\beta)$ qui sont tous deux
des isomorphismes. Le compl\'ementaire de $E$ dans~$\overline E$ est
alors constitu\'e de l'{\og infini\fg} (donn\'e par $\alpha=0$) et de {\og
z\'ero\fg}
(donn\'ee par $\beta=0$).


Ainsi, les projections de 
$\mathscr W$ vers $ \mathscr O_{A^0}$ (resp.\ 
$\mcL^\vee)$ d\'efinissent deux sous-sch\'emas de $\overline E$, respectivement
les sections {\og nulle\fg} et {\og infini\fg} (la
section nulle est effectivement la section nulle de $\V(\mcL^\vee)$).
Ce sont deux diviseurs relatifs de $\overline E$ au-dessus de $A^0$, 
notons les $D_0$ (resp. $D_\infty$).
Notant $\pi$ la projection $\P(\mathscr W)\ra A$,
il r\'esulte du lemme suivant que 
$\mathscr O_{\P(\mathscr W)}(D_0-D_\infty)=\pi^*\mcL$.

\begin{lemm} \label{l:fibre}
(cf.~\cite[Chap. V, Prop. 2.6]{hartshorne77})\quad
Soient $X$ un sch\'ema, $\mathscr E$ un faisceau localement libre
de rang $n+1$ sur $X$
et $\pi:\P=\P(\mathscr E)\ra X$. Si $\mathscr N$ et $\mathscr V$
sont deux faisceaux localement libres sur $X$,
de rang $1$ et $n$ respectivement,
avec une suite exacte $0\ra\mathscr N \ra \mathscr E\ra \mathscr V\ra 0$,
l'image de $\P(\mathscr V)\hookrightarrow \P$ est un diviseur $D$ dans $\P$
et $\mathscr O_\P(1) =\mathscr O_\P(D)\otimes \pi^*\mathscr N$.
\end{lemm}
\begin{proof} 
Posons $\mathscr E'=\mathscr E\otimes \mathscr N^\vee$ et $\P'=\P(\mathscr E')$.
Comme $\mathscr N$ est inversible, $\P'$ est canoniquement isomorphe \`a $\P$,
le faisceau $\mathscr O_{\P'}(1)$ s'identifiant 
d'apr\`es~\cite[Chap. II, Lemma 7.9]{hartshorne77} \`a 
$\mathscr O_\P(1)\otimes\mathscr N^{\vee}$.
Cela nous ram\`ene \`a prouver le lemme quand
$\mathscr N$ est trivial. Dans ce cas, l'injection $\mathscr O_X\hra \mathscr E$
s'interpr\`ete comme un \'el\'ement non nul de $\Gamma(X,\mathscr E)$, puis
comme $\pi_*\mathscr O(1)=\mathscr E$, comme une section $s$ non nulle de 
$\Gamma(\P,\mathscr O_\P(1))$ dont le diviseur est \'egal \`a $D$~; ainsi,
le lemme est d\'emontr\'e.
\end{proof}

Notons $\mathscr M_0$ et $\mathscr M_\infty$ les faisceaux inversibles
associ\'es aux diviseurs $D_0$ et $D_\infty$. Ainsi, $\mathscr M_0\otimes
\mathscr M_\infty^\vee \simeq \pi^*\mcL$.
Si $\sigma:\Spec \C\ra S$ est un point complexe de $S$, montrons comment
munir les faisceaux inversibles $\sigma^*\mathscr M_0$ (resp. 
$\sigma^*\mathscr M_\infty$) de m\'etriques hermitiennes.

Montrons tout d'abord l'existence d'une {\og m\'etrique carr\'ee\fg}
sur $\mcL$ (cf.~\cite[II.2]{moret-bailly85b} dans le cas cubiste)~:
\begin{prop}
Soient $A$ une vari\'et\'e ab\'elienne complexe et $\mcL$ un faisceau
inversible sur $A$ algébriquement équivalent à zéro
et rigidifi\'e \`a l'origine.
Alors, il existe une unique m\'etrique hermitienne sur $\mcL$ telle que
l'unique  {\og isomorphisme carr\'e\fg}
$m^*\mcL\simeq p_1^*\mcL\otimes p_2^*\mcL$
compatible avec la rigidification à l'origine soit une isométrie.
Cette m\'etrique est de plus
l'unique m\'etrique hermitienne sur $\mcL$ compatible \`a la rigidification
et dont la forme de courbure est nulle.
\end{prop}
\begin{proof}
Tout d'abord, $c_1\in H^2_{dR}(A)$
est la premi\`ere classe de Chern de $\mcL$, il existe d'apr\`es la th\'eorie
de Hodge une unique forme diff\'erentielle invariante par
translations  qui repr\'esente $c_1$.
D'autre part, le {\og lemme $\partial\overline\partial$\fg}
(cf.~\cite[pp.\ 148--149]{griffiths-h78}) implique l'existence d'une
m\'etrique hermitienne sur $\mcL$ dont la forme de courbure soit cette
forme diff\'erentielle, et deux telles m\'etriques diff\`erent d'une
constante strictement positive. Il existe ainsi sur $\mcL$ 
une unique m\'etrique qui soit compatible \`a la trivialisation \`a l'origine
et dont la forme de courbure soit invariante par translations.
D'autre part, $\mcL$ appartenant \`a
$\Pic^0(A)$, on a $c_1=0$ et la courbure de la m\'etrique est nulle.

Enfin, le fibr\'e $m^*\mcL\otimes p_1^*\mcL^\vee\otimes p_2^*\mcL^\vee$,
trivial, est muni d'une m\'etrique hermitienne dont la forme de courbure
est nulle. Par suite, il poss\`ede une section globale sans z\'eros
dont la norme est
une fonction harmonique et donc constante, $A$ \'etant compacte. Ainsi,
la norme de l'isomorphisme carr\'e est constante~; sa valeur \`a l'origine
est par d\'efinition \'egale \`a~$1$, d'o\`u la proposition.
\end{proof}

\begin{coro}\label{coro:isom}
Avec les notations de la proposition pr\'ec\'edente, l'unique isomorphisme
$[n]^*\mcL \simeq \mcL^{\otimes n}$ compatible aux rigidifications \`a
l'origine est une isom\'etrie.
\end{coro}
\begin{proof}
La proposition pr\'ec\'edente nous fournit sur le faisceau inversible
$[n]^*\mcL\otimes \mcL^{\vee n}$,
canoniquement isomorphe au faisceau $\mathscr O_A$,
une m\'etrique hermitienne canonique
dont il faut v\'erifier qu'elle est triviale.
Or, d'une part cette m\'etrique est constante (sa forme de courbure \'etant
nulle), et d'autre part, la norme de la section 1 vaut 1 \`a l'origine,
ce qui ach\`eve la preuve du corollaire.
\end{proof}

La proposition pr\'ec\'edente nous fournit une m\'etrique canonique sur
$\sigma^*\mcL$, si bien que $\mathscr W$ est muni,
pour tout point complexe $\sigma:\Spec\C\ra S$ de $S$, d'une m\'etrique
continue $\|\cdot\|_{\mathscr W}^\sigma$~: 
si $s=(s_1,s_2)$ est une section locale
de $\mathscr O_A\oplus\mcL^\vee$, on d\'efinit
$$ \|s\|^\sigma_{\mathscr W} (x) 
=  \|s_1\|(x^\sigma) + \|s_2\|(x^\sigma) . $$

Comme le faisceau inversible
$\mathscr O_\P(1)$ est un quotient de $\pi^*\mathscr W$, il est
naturellement muni d'une m\'etrique hermitienne.
D'apr\`es le lemme~\ref{l:fibre}, $\mathscr M_0 = \pi^*\mcL\otimes\mathscr
O_\P(1)$ ce qui nous donne une m\'etrique hermitienne canonique $\omega_0$
sur $\mathscr M_0$ en prenant le produit tensoriel des
deux m\'etriques sur $\pi^*\mcL$ et sur $\mathscr O_\P(1)$. 
De m\^eme, $\mathscr M_\infty$ est muni d'une
m\'etrique hermitienne canonique $\omega_\infty$.

Donnons maintenant une formule explicite pour la norme
en un point $P\in E(\C)$
des sections canoniques $s_{D_\infty}$ et $s_{D_0}$ des faisceaux
$\mathscr M_\infty$ et $\mathscr M_0$
dont le diviseur est $D_\infty$ et $D_0$.

\begin{lemm}\label{lemm:calcul}
Fixons une place complexe $\sigma:\Spec\C\ra S$.
Soit $x\in A(\C)$
et $e$ une base norm\'ee de $\mathscr L_x^\vee$~; soit
aussi un point $P\in\overline E(\C)$ relevant $x$, ainsi
qu'une base $\eps$ de $\mathscr O_\P(1)_P$.
Le point $P$ correspond alors \`a deux nombres complexes $u_1$ et $u_2$
par le quotient
$$\mathscr O_x \oplus \mathscr L_x^\vee \rightarrow \mathscr O_{\P}(1)_P,
\qquad
    t_1+t_2 e \mapsto (t_1 u_1 + t_2 u_2 ) \eps.
$$
Alors,
$$ \|{s_{D_\infty}}\|(P) = \frac{|u_1|}{\max(|u_1|,|u_2|)}
\quad\text{et}\quad
\| s_{D_0}\|(P) = \frac{|u_2|}{\max(|u_1|,|u_2|)}. 
$$
\end{lemm}
\begin{rema}
Ces m\'etriques sont seulement continues alors que la g\'eom\'etrie
d'Arakelov consid\`ere usuellement des m\'etriques $\mathscr C^\infty$ ; 
c'est cependant cette m\'etrique qui refl\`ete pr\'ecis\'ement 
l'action du tore sur la compactification, cf.\ la
proposition~\ref{prop:n-iso}, et donnera ainsi lieu aux
hauteurs canoniques.
D'autre part, lorsqu'on consid\`ere la hauteur de points rationnels,
il suffit de choisir une m\'etrique continue.
Enfin, comme c'est une limite uniforme de m\'etriques lisses, les arguments
de~\cite{zhang95} montrent que la consid\'eration de cette m\'etrique
est l\'egitime dans le contexte de la g\'eom\'etrie d'Arakelov en dimension
sup\'erieure, par exemple pour \'etudier la hauteur des cycles.
\end{rema}
\begin{proof}[Preuve du lemme]
La section
$(1,0)\in\Gamma(\pi^*(\mathscr O_{A}\oplus\mathscr L^\vee))$
ayant pour image la section $s_{D_\infty}$ de $\mathscr O_\P(1)$,
on a
\begin{eqnarray*}
\|{s_{D_\infty}}\|(P)
 &=& \inf_{(x_1,x_2 e)\mapsto s_\infty} \left\|x_1+x_2 e\right\| _{\mathscr
W}
 = \inf _{x_1u_1+x_2u_2 = u_1} |x_1|+|x_2| \\
 &=& |u_1| \inf_{(x_1,x_2)\in\C^2}  \frac{|x_1|+|x_2|}{|x_1 u_1+x_2 u_2|}
 =  \frac{|u_1|}{\max(|u_1|,|u_2|)}.
\end{eqnarray*}

La section rationnelle $s_{D_0}/s_{D_\infty}$ de $\mathscr
M_0\otimes\mathscr M_\infty^\vee=\pi^*\mathscr L$ associe
au point $P$ l'\'el\'ement $(u_2/u_1)e^\vee$. Sa norme est donc
$|u_1/u_2|$ puisque $\norm{e^\vee}=\norm{e}=1$. On a donc
$$ \norm{s_{D_0}}(P) = |u_2/u_1| \norm{s_{D_\infty}}(P)
= \frac{|u_2|}{\max(|u_1|,|u_2|)}. \qed $$
\let\qedsymbol\relax
\end{proof}

\'Etudions enfin le comportement des objets que nous venons d'introduire
par rapport aux morphismes de multiplication par $n$ dans $E$.
\begin{prop}\label{prop:n-iso}
Le morphisme $[n]_E:E\ra E$ s'\'etend uniquement en un morphisme $\overline E\ra
\overline E$, toujours not\'e $[n]$. De plus, 
on a des isomorphismes canoniques de norme~$1$~: si $n\geq 0$,
$ [n]^*\mathscr M_0 \simeq \mathscr M_0^n$ et
$ [n]^*\mathscr M_\infty \simeq \mathscr M_\infty^n$~;
si $n\leq 0$, on a en revanche 
$[n]^*\mathscr M_0 \simeq \mathscr M_\infty^{|n|}$
et $[n]^*\mathscr M_\infty \simeq \mathscr M_0^{|n|}$.
\end{prop}
\begin{proof}
Si $n\in\Z$, la multiplication par $n$ dans $E$ provient du diagramme
suivant~:
$$
\xymatrix{
  {\V(\mcL^\vee)\setminus\{0\}} \ar[r] \ar[rrd] 
  & {\V(\mcL^{\vee n})\setminus\{0\}} \ar@{=}[r] \ar[rd] 
  & {\V([n]^*\mcL^\vee)\setminus\{0\}} \ar[r] \ar[d] \ar@{} [rd] |{\square}
  & {\V(\mcL^\vee)\setminus\{0\}} \ar[d] \\
  & & {A^0} \ar[r] & {A^0\rlap{\quad,}} }
$$
o\`u, 
l'application $\V(\mcL^\vee)\setminus\{0\}\ra\V(\mcL^{\vee n})\setminus\{0\}$ 
associe \`a une section sans z\'eros de $\mcL$ la puissance tensorielle 
$n$-\`eme\ de cette section, et le carr\'e de droite est cart\'esien.
Il en r\'esulte que le morphisme $[n]:E\ra E$ s'\'etend \`a $\overline E$
selon le diagramme, dont le carr\'e de droite est cart\'esien~:
$$
\xymatrix{
  {\P(\mathscr W)} \ar[r] \ar[rrd] 
  & {\P(\mathscr O_{A^0}\oplus\mcL^{\vee n})} \ar@{=}[r] \ar[rd]
  & {\P([n]^*\mathscr W)} \ar[r] \ar[d] \ar@{} [rd] |{\square}
  & {\P(\mathscr W)} \ar[d] \\
  & & {A^0} \ar[r] & {A^0\rlap{\quad,}} }
$$
la fl\`eche $\P(\mathscr W)\ra  \P(\mathscr O_{A^0}\oplus\mcL^{\vee n})$
\'etant donn\'ee 
au niveau des~$S$-points par l'application
$$ \big( (\alpha,\beta):\mathscr O_S\oplus \eps_P^*\mcL^\vee\ra\mathscr J
\big)
  \longmapsto
   \big( (\alpha^{\otimes n},\beta^{\otimes n}):\mathscr
O_S\oplus \eps_P^*\mcL^{\vee n} \ra\mathscr J^{\otimes n} \big)
$$
quand $n\geq 0$, et par
\begin{multline*}
\big( (\alpha,\beta):\mathscr O_S\oplus \eps_P^*\mcL^\vee\ra\mathscr J \big)
  \longmapsto \\
   \big( (\beta^{\otimes |n|},\alpha^{\otimes |n|})\otimes
\id_{\mcL^{\otimes |n|}}:
\mathscr O_S\oplus \eps_P^*\mcL^{\otimes |n|} \ra
\mathscr J^{\otimes |n|}\otimes\mcL^{\otimes|n|} \big) 
\end{multline*}
lorsque $n\leq 0$.

En g\'en\'eral, un morphisme $\P(\mathscr W)\ra\P(\mathscr W)$ relevant
la multiplication par $n$ sur $A^0$ qui envoie par image r\'eciproque le
faisceau $\mathscr O_\P(1)$ sur le faisceau $\mathscr O_\P(n)$
correspond \`a la donn\'ee d'un morphisme surjectif
$$ \pi^* [n]_A^* (\mathscr O_{A^0}\oplus\mcL^\vee) \ra \mathscr O_\P(n). $$
La multiplication par $n\geq 0$ sur $\overline E$ est ainsi donn\'ee par
les fl\`eches naturelles
$$ \pi^* (\mathscr O_{A^0} \oplus  \mcL^{\vee n} ) 
\ra \pi^*\Sym^n (\mathscr O_{A^0} \oplus \mcL^\vee ) \ra \mathscr O_\P(n). $$
Lorsque $n\leq 0$, la multiplication par $n$ correspond \`a la composition
\begin{multline*}
\pi^*(\mathscr O_{A^0}\oplus \mcL^{\otimes |n|} )
\simeq \pi^* (\mcL^{\vee n}\oplus \mathscr O_{A^0}) \otimes\pi^*
\mcL^{\otimes |n|} \ra \\
\ra \pi^* \Sym^{|n|} (\mcL^\vee\oplus\mathscr
O_{A^0})\otimes\pi^*\mcL^{\otimes |n|}
\ra \mathscr O_\P(|n|)\otimes\pi^*\mcL^{\otimes|n|}.
\end{multline*}

Quand $n\geq 0$, on a ainsi $[n]_{\overline E}^*\mathscr O_\P(1)=\mathscr
O_\P(n)$, tandis que quand $n\leq 0$, on a
$ [n]_{\overline E}^*\mathscr O_\P(1)=\mathscr
O_\P(|n|)\otimes\pi^*\mcL^{\otimes|n|}$.
 

Il reste \`a montrer que ces isomorphismes respectent les m\'etriques
hermitiennes~: pour cela, il suffit de montrer que les isomorphismes
$$
[n]^*\mathscr O_\P(1) = \mathscr O_\P(n) \quad\text{(pour $n\geq 0$)} 
\qquad\text{et}\qquad 
[-1]^*\mathscr O_\P(1)=\mathscr O_\P(1)\otimes \pi^*\mcL
$$
respectent eux-m\^emes les m\'etriques, les formules pour $\mathscr M_0$
s'en d\'eduiront puisque $[n]^*\mathscr L$ est isom\'etrique \`a $\mathscr
L^n$ (corollaire~\ref{coro:isom}).

Dans un souci d'all\`egement,
on effectue le changement de base de $S$ \`a $\C$ sans changer
les notations.
Soient $x$ un point de $A(\C)$ et $P\in\bar E(\C)$ relevant $x$, et,
comme dans le lemme~\ref{lemm:calcul}, $e$ une base norm\'ee
de $\mathscr L^\vee_x$, $\eps$ une base de $\mathscr O(1)_P$,
et $(u_1,u_2)\in\C^2$ tels que $P$ soit d\'efini par le quotient
$\mathscr O_x\oplus \mathscr L_x^\vee\ra \mathscr O_P(1)_P$,
$x_1+x_2 e\mapsto (x_1u_1+x_2u_2)\eps$.

Alors, si $n\geq 1$,
$f=e^{\otimes n}$ s'identifie \`a un \'el\'ement non nul
de $\mathscr L_{[n]x}^\vee$, dont la norme est $\|f\|=\| e\| ^n=1$.
De plus, dans les bases $f$ de $\mathscr L^\vee_{[n]x}$
et $\eps^{\otimes n}$ de $\mathscr O_\P(1)_{[n]P}$,
le point $[n]P$ correspond au couple $(u_1^n,u_2^n)$ si bien
que
$ \|s_\infty\| ([n]P) = \| s_\infty\| (P)^n$,
ainsi qu'il fallait d\'emontrer.

Pour $n=-1$, soit $f=e^{\vee}$ la base de $\mathscr L_{-x}^\vee\simeq
\mathscr L_x^{\vee\vee}$ duale de $e$, de sorte que $\norm{f}=1$.
Le point $[-1]P$ correspond au quotient
$\mathscr O_x\oplus \mathscr L_{-x}^\vee \ra \mathscr
O_P(1)_P\otimes\mathscr L_x$ d\'efini par
$(x_1+x_2 f) \mapsto (x_1 u_2 + x_2 u_1)\eps\otimes f$, de sorte que
$$ \| s_{D_\infty}\| ([-1]P)
= \frac{ |u_2| }{\max(|u_1|,|u_2|)}
= \| s_{D_0} \| (P) 
$$
d'o\`u le r\'esultat puisque $[-1]^* D_0=D_\infty$.
\end{proof}

%
%

\begin{rema}\label{rema:semiab}
La m\^eme m\'ethode permet de traiter le cas d'une vari\'et\'e
semi-ab\'elienne dont le tore sous-jacent est d\'eploy\'e. En
effet, cela nous ram\`ene \`a une extension d'une vari\'et\'e
ab\'elienne par une puissance $\gm^t$, d'o\`u $t$ fibr\'es
alg\'ebriquement \'equivalents \`a~$0$~: $\mcL_1$, \dots, $\mcL_t$
que l'on peut m\'etriser comme pr\'ec\'edemment.
On dispose alors (entre autres)
de deux compactifications naturelles, \`a savoir
$\P(\mathscr O_A\oplus \mcL_1^\vee)\times_A \cdots \times_A \P(\mathscr
O_A\oplus\mcL_t^\vee)$
et $\P(\mathscr O_A\oplus \mcL_1^\vee\oplus \cdots \oplus\mcL_t^\vee)$.
Dans l'un et l'autre cas, on dispose de faisceaux inversibles
m\'etris\'es construits \`a partir des $\mathscr O(1)$
et des $\mcL_i$. Ils donneraient lieu \`a des hauteurs canoniques,
comme au paragraphe suivant.
\end{rema}

\section{Construction des hauteurs relatives}
\label{sec:t-haut}
On reprend les notations du paragraphe pr\'ec\'edent, en supposant
que $S$ est le spectre de l'anneau des entiers d'un corps de nombres~$K$
Rappelons que le premier groupe de Chow--Arakelov $\hPic(S)$ de $S$
s'identifie au groupe des classes d'isomorphisme
de faisceaux inversibles sur $S$ munis de
m\'etriques hermitiennes {\og \`a l'infini\fg} compatibles \`a
la conjugaison complexe.

Fixons tout d'abord
un entier $N>0$ qui annule les groupes des composantes connexes
de $A_s$ pour tout point $s\in S$.

Soient $P\in E(\eta)$ et $Q=\pi(P)\in A(\eta)$.
Comme $A/S$ est le mod\`ele de N\'eron de $A_\eta$, il existe une section
$\eps_Q:S\ra A$ qui prolonge $Q$. 
D'apr\`es le choix de l'entier $N$, le point
$[N]_A Q\in A(\eta)$ se prolonge en une
section $\eps_{[N]Q}:S\ra A^0$. Par suite, $\overline E/A^0$ \'etant projectif,
$[N]P$ se prolonge en une section $\eps_{[N]P}:S\ra \overline E$ qui rel\`eve
$\eps_{[N]Q}$.

\begin{prop}\label{prop:n-iso2}
Les \'el\'ements $(\eps_{[N]P}^* \mathscr M_0)\otimes \frac1N$ 
et  $(\eps_{[N]P}^* \mathscr M_\infty)\otimes \frac1N$ de
$\hPic(S)\otimes_\Z \Q$ ne d\'ependent pas du choix de $N$. 
On les note respectivement $H_{0,S}(P)$ et $H_{\infty,S}(P)$.
\end{prop}
\begin{proof}
Fixons $\sbullet$ l'un des symboles $\{0,\infty\}$.
Le caract\`ere canonique des m\'etriques hermitiennes $\omega_\sbullet$
sur les faisceaux
inversibles $\mathscr M_\sbullet$ implique qu'elles sont invariantes par
la conjugaison complexe. Ainsi, nous avons bien par fonctorialit\'e
des \'el\'ements $\eps_{[N]P}^*\mathscr M_\sbullet$ dans $\hPic(S)$.

D'autre part, si $M$ est un autre entier qui annule les groupes des
composantes connexes de $A_s$ pour tout $s\in S$, montrons que
$$(\eps_{[N]P}^*\mathscr M_\sbullet)\otimes \frac1N 
	= (\eps_{[M]P}^*\mathscr M_\sbullet ) \otimes \frac 1M. $$
Pour cela, on peut supposer que $M$ est un multiple de $N$, soit
$M=Nk$ pour un entier $k\geq 1$.
Or d'une part, $[k]^*\mathscr M_\sbullet=\mathscr M_\sbullet ^{\otimes k}$ en tant
que faisceau inversible m\'etris\'e (proposition~\ref{prop:n-iso})
et d'autre part,
$\eps_{[M]P}=[k]\circ \eps_{[N]P}$, si bien que l'on a 
$$\eps_{[M]P}^*\mathscr M_\sbullet  = (\eps_{[N]P}^*\mathscr M_\sbullet)^{\otimes k}, $$
ce qui conclut la preuve de la proposition.
\end{proof}

D'autre part, la proposition~\ref{prop:n-iso} (ou la 
proposition~\ref{prop:n-iso2}, comme on veut~!) entra\^\i ne imm\'ediatement
la proposition suivante~:
\begin{prop}\label{prop:n-can}
Soit $P\in E(\eta)$. Si $n\geq 0$, on a $H_{0,S}([n]P)=n 
H_{0,S}(P)$ et $H_{\infty,S}([n]P)=nH_{\infty,S}(P)$.
De plus $H_{0,S}([-1]P)=H_{\infty,S}(P)$.
\end{prop}

\begin{prop}
Soient $\eta'\ra \eta$ une extension finie, $f:S'\ra S$ le normalis\'e de 
$S$ dans $\eta'$, $A'$ le mod\`ele de N\'eron de $A_\eta\times_\eta\eta'$,
$E'$ l'extension de $A^{\prime 0}$ par $\gm$ qui prolonge
$E_\eta\times\eta'$.
Si $P\in E_\eta(\eta)$, on a 
$$ H_{\infty,S'}(P\times_S S') = f^* H_{\infty,S}(P) 
\in \hPic(S') \otimes_\Z \Q \quad , $$
et de m\^eme pour $H_{0}$.
\end{prop}
\begin{proof}
Si $A/S$ est semi-stable, c'est clair~: la composante neutre de $A'$
est obtenue \`a partir de celle de $A$ par changement de base, si bien
que $E'=E\times_S S'$, etc.

Dans le cas g\'en\'eral, soit $\phi:A\times_S S'\ra A'$ le morphisme 
naturel qui prolonge l'identit\'e $A_\eta\times \eta'=A_{\eta'}$. 
L'image de $A^0\times_S S'$ par $\phi$ est contenue dans 
$A^{\prime 0}$  et il nous faut comparer $\phi^*\mcL'$ et $\mcL\times_S S'$.
Or, $\phi^*\mcL'\otimes f^*\mcL^\vee$ est un faisceau inversible 
m\'etris\'e sur $A^0\times_S S'$ qui v\'erifie le th\'eor\`eme du carr\'e
et est g\'en\'eriquement trivial. Comme $S'\ra S$ est fid\`element plat,
on a $(\pi\times_S S')_*\mathscr O_{A^0\times_S S'}=\mathscr O_{S'}$~;
d'autre part, $A^0\times_S S'\ra S'$ est \`a fibres connexes, si bien que
$\phi^*\mcL'\otimes f^*\mcL^\vee$ provient d'un faisceau inversible sur $S'$,
lequel est trivial \`a cause des rigidifications.
Autrement dit, $\phi^* E'$, etc.
sont obtenues \`a partir de $E$ par changement de base, d'o\`u la 
proposition.
\end{proof}

Nous pouvons donc poser~:
\begin{defi}
Soient $\overline\eta$ la cl\^oture alg\'ebrique de $\eta$ et $P\in E(\overline\eta)$.
Si $\eta'$ est une extension finie de $\eta$ telle que $P\in E(\eta')$
et $f:S'\ra S$ est le normalis\'e de $S$, on appelle {\em hauteurs 
relatives} de $P$ les r\'eels
$\hdeg H_0(P):=\frac{1}{[S':S]} \hdeg H_{0,S'}(P)$
et $\hdeg H_\infty(P):= \frac{1}{[S':S]}\hdeg H_{\infty,S'}(P)$.
\end{defi}

\begin{theo}\label{theo:t-hcan}
Les fonctions $\hdeg H_0$ et $\hdeg H_\infty : E(\overline\eta)\ra \R$ 
sont les hauteurs canoniques sur $\overline E(\overline\eta)$
attach\'ees aux faisceaux inversibles $\mathscr M_{0,\eta}$ et $\mathscr
M_{\infty,\eta}$ sur $\overline E_\eta$~; elles sont positives.
De plus (Zarhin--Bloch, Mazur--Tate), 
$\hdeg H_0(P) - \hdeg H_\infty(P)$
est la hauteur de N\'eron--Tate de $\pi(P)$ relativement au faisceau
alg\'ebriquement \'equivalent \`a z\'ero $\mcL$ sur $A$.
\end{theo}
\begin{proof}
Pour $\sbullet\in \{0,\infty\}$,
soit $h_\sbullet$ une hauteur de Weil sur $\overline E(\overline\eta)$
attach\'ee
au faisceau inversible $\mathscr M_{\sbullet,\eta}$ sur $\overline E_\eta$.
Fixons un entier $n\geq 2$.
Les {\em hauteurs canoniques} relatives au faisceau inversible
$\mathscr M_{\sbullet,\eta}$ de poids $n$ pour le morphisme $[n]_{\overline E_\eta}$
(cf.\ par exemple~\cite{call-s93})
sont par d\'efinition les fonctions sur
$\overline E(\overline\eta)$ d\'efinies par 
$$ \hat h_\sbullet (P)= \lim_{k\ra\infty} \frac{1}{n^k} h_\sbullet([n^k]P). $$
Rappelons pourquoi cette limite existe~: comme $[n]^*\mathscr 
M_{\sbullet,\eta}=\mathscr M_{\sbullet,\eta}^n$, il existe une constante $C_n$
telle que $|h_\sbullet([n]P)-n h_\sbullet(P)|\leq C_n$ pour tout point $P$,
si bien que
$$\left| 
\frac{1}{n^k} h_\sbullet([n^k]P)- \frac{1}{n^{k-1}} h_\sbullet([n^{k-1}]P) 
  \right| \leq \frac{1}{n^k} C_n, $$
et la limite s'\'ecrit comme somme d'une s\'erie uniform\'ement 
convergente. La d\'enomination hauteur est justifi\'ee par la comparaison
$|\hat h_\sbullet- h_\sbullet|\leq C_n/(n-1)$, tandis que le terme canonique
vient de ce que $\hat h_\sbullet([n]P)=n\hat h_\sbullet(P)$ pour tout $P$.

Montrons maintenant qu'il existe pour toute extension finie $\eta'\ra\eta$
une constante $C_{\eta'}$ telle que l'on ait, pour
tout point $P\in \overline E(\eta')$, l'in\'egalit\'e
$$ \left| h_\sbullet (P) - \hdeg H_\sbullet (P) \right| \leq C_{\eta'}. $$
En effet, si $S'$ est le normalis\'e de $S$ dans $\eta'$
et si l'entier $N>0$ annule les groupes des composantes connexes du mod\`ele
de N\'eron de $A_{\eta'}$ sur $S'$, 
l'expression $h_\sbullet ([N]P)-N h_\sbullet(P)$
est born\'ee uniform\'ement en $P\in E(\eta')$.
D'autre part, $\hdeg H_\sbullet([N]P)=N \hdeg H_\sbullet(P)$
d'apr\`es la proposition~\ref{prop:n-can}. Enfin, en choisissant comme
mod\`ele entier de $\overline E_{\eta'}$ l'adh\'erence de $E'$ dans un espace
projectif convenable, on constate que la diff\'erence
$\hdeg H_\sbullet([N]P)-h_\sbullet([N]P)$ est uniform\'ement born\'ee
lorsque $P$ d\'ecrit $\overline E(\eta')$.
En mettant bout \`a bout ces majorations, on a bien une in\'egalit\'e
comme annonc\'ee.

La proposition~\ref{prop:n-can} et la d\'efinition de la hauteur canonique
entra\^{\i}nent alors que pour tout $P\in \overline E(\eta')$,
$$ \hat h_\sbullet (P) = \hdeg H_\sbullet (P). $$
L'extension finie $\eta'$ \'etant arbitraire, cela 
implique bien que $\hat h_\sbullet = \hdeg H_\sbullet$.

Comme $\mathscr M_\infty$ est effectif,
la hauteur $h_{\infty}$ est minor\'ee sur le compl\'ementaire
de $D_\infty$, ce qui implique que $\hat h_{\infty}=\hdeg H_\infty$ est
positive. De m\^eme pour $\hdeg H_0$.
(Il est aussi possible d'utiliser le fait que les sections $s_{D_0}$
et $s_{D_\infty}$ de $\mathscr M_0$ et $\mathscr M_\infty$
ont une norme $\leq 1$ en tout point~; voir les formules
du \S 6.)

Enfin, comme $\mathscr M_0\otimes \mathscr M_\infty^\vee=\pi^*\mcL$,
on a pour tout $P\in \overline E(\eta)$ relevant un point de $A^0$,
$$ \hdeg H_0(P) - \hdeg H_\infty(P) = \hdeg \eps_P^*\pi^*\mcL
 = \hdeg (\pi\circ\eps_P)^*\mcL. $$
Autrement dit,  $h_\mcL$ d\'esignant la hauteur de N\'eron--Tate
sur $A_\eta$ relative au fibr\'e inversible $\mcL$,
la fonction lin\'eaire
$$ \hdeg H_0-\hdeg H_\infty- h_\mcL\circ \pi$$
est, pour tout $\eta'\ra\eta$, born\'ee sur un sous-groupe d'indice fini
de $A(\eta')$~; elle est alors n\'ecessairement nulle, ce qu'il fallait
d\'emontrer.
\end{proof}

Remarquons pour finir que les hauteurs relatives contiennent toute
l'information n\'ecessaire pour conna\^\i tre la hauteur d'un point
de l'extension (compactifi\'ee)
dans un plongement projectif donn\'e~: le groupe
de Picard de $\overline E_\eta$ est $\Z\oplus \Pic(A_\eta)$, si bien que tout
faisceau (tr\`es) ample sur $\overline E_\eta$ est de la forme $\pi^*\mathscr N
\otimes \mathscr O_\P(n)$ pour un faisceau (tr\`es) ample $\mathscr N$
sur $A_\eta$ et un entier $n>0$.
En sus de la hauteur de N\'eron--Tate sur $A_\eta$,
la connaissance de $H_\infty$ suffit donc, m\^eme si la consid\'eration
de $H_0 + H_\infty$ est plus sym\'etrique.

\begin{coro}[{Waldschmidt, \protect\cite[App.]{cohen86}, 
  Laurent~\protect\cite{laurent80}}]
\label{coro:t-hcan}
Soit $h_{NT}$ une hauteur de N\'eron--Tate sur $A_\eta$ pour un diviseur
sym\'etrique ample et $h$ la fonction $E_\eta(\overline\eta)\ra\R_+$ d\'efinie par
$h_{NT}\circ\pi + \hdeg H_0 + \hdeg H_\infty$. C'est une hauteur sur
$E_\eta$~; de plus $h(P)=0$ si et seulement si $P$ est d'ordre fini.
\end{coro}
\begin{proof}
Que ce soit une hauteur r\'esulte du th\'eor\`eme~\ref{theo:t-hcan}
et des remarques qui pr\'ec\`edent~;
elle est positive comme somme de fonctions positives.
Enfin, si $h(P)=0$, il est n\'ecessaire que $h_{NT}(\pi(P))=0$,
ce qui prouve que $\pi(P)$ est d'ordre fini. Alors, pour tout entier
$n$, $h([n]P)=nh(P)=0$ et le th\'eor\`eme de Northcott
entra\^\i ne que l'ensemble $\{P,2P,3P,\ldots\}$ est fini, c'est-\`a-dire
que $P$ est d'ordre fini. La r\'eciproque est claire, d'o\`u le corollaire.
\end{proof}

\section{Points de hauteur relative nulle}
\label{sec:t-misc}
Pour finir, nous voulons donner quelques expressions explicites des
hauteurs $H_0$ et $H_\infty$ et appliquer la th\'eorie pr\'ec\'edente
\`a l'\'etude des points de hauteur relative nulle. 

Nous conservons les notations des paragraphes pr\'ec\'edents.
Rappelons qu'un point $P$ de $\overline E=\P(\mathscr W)$
relevant une section (de la composante neutre) $\eps_Q:S\ra A^0$
est la donn\'ee d'un quotient inversible de rang~$1$ de $\eps_Q^*\mathscr W$,
soit un faisceau inversible $\mathscr I$ sur $S$ et deux sections
$x_0:\mathscr O_S\ra \mathscr I$, $x_1:\mathscr \eps_Q^*\mcL^\vee\ra \mathscr
I$ telles que $(x_0,x_1):\eps_Q^*\mathscr W\ra\mathscr I$ soit surjectif.
Le point $P_\eta$ appartient \`a~$E_\eta$ si et seulement
si $x_0\neq 0$ et $x_1\neq 0$.

Les faisceaux inversibles $\mathscr M_\infty$ et $\mathscr M_0$ \'etant
d\'efinis par les inclusions $\mathscr O_{A^0}\hra \mathscr W$ et $\mcL^\vee\hra
\mathscr W$, les diviseurs sur $S$ des \'el\'ements  
$\theta_\infty=x_0\in\Hom(\mathscr O_S,\mathscr I)$
et $\theta_0=x_1\in\Hom(\eps_Q^*\mcL^\vee,\mathscr I)$ 
s'interpr\`etent respectivement comme les intersections (propres,
sur le sch\'ema $\overline E$ puisqu'on a suppos\'e que $P_\eta$
\'etait un point de $E_\eta$)
$D_\infty\cdot \overline{\{P\}}$ et $D_0\cdot\overline{\{P\}}$.

Alors, $\eps_P^*\mathscr M_\infty = \mathscr I$
et la norme de la section $\theta_\infty$ de $\mathscr I$ pour la
m\'etrique hermitienne h\'erit\'ee de $\mathscr M_\infty$,
en tout point $\sigma\in S(\C)$, est donn\'ee par la formule
$$ \| \theta_\infty\|^\sigma 
= \frac{\|x_0\|^\sigma}{\max(\|x_0\|^\sigma,\|x_1\|^\sigma)} $$
(\`a proprement parler, le membre de droite pour \^etre calcul\'e,
n\'ecessite le choix d'une m\'etrique sur $\mathscr I$ mais n'en d\'epend pas).
Notons $\div(\theta_\infty)
=\sum\limits_{\fp \in S} n_{\infty,\fp }[\fp ]$,
on a alors
$$ \hdiv(\theta_\infty)= \left(\div(\theta_\infty),
-\log\|\theta_\infty\|^{2,\sigma} \right) $$
et 
$$
\hdeg H_\infty(P) = \hdeg\,\hdiv (\theta_\infty) = \sum_{\fp }
n_{\infty,\fp }\log N(\fp ) - \sum_{\sigma\in S(\C)}
\log\|\theta_\infty\|^\sigma. 
$$
C'est une somme de termes positifs de sorte que la hauteur relative $\hdeg
H_\infty(P)$ est nulle si et seulement si $\theta_\infty$ est un
isomorphisme dans $\hPic(S)$~: $\theta_\infty$ est un isomorphisme de
faisceaux inversibles et sa norme est 1 en toute place~; cette derni\`ere
condition \'equivaut \`a $\|x_0/x_1\|^\sigma\geq 1$.

De m\^eme, si $\div(\theta_0)=\sum_{\fp } n_{0,\fp } [\fp ]$, on a
$$
\hdeg H_0(P)=\sum_{\fp \in S} n_{0,\fp } \log N(\fp ) -
\sum_{\sigma\in S(\C)} \log \|\theta_0\|^\sigma  \  , 
$$
o\`u
$$ \|\theta_0\|^\sigma 
= \frac{\|x_1\|^\sigma}{\max(\|x_0\|^\sigma,\|x_1\|^\sigma)}. $$
Ainsi, la hauteur relative $\hdeg H_0(P)$ est nulle si et seulement si 
$\theta_0$ est un isomorphisme dans $\hPic(S)$, soit si c'est un
isomorphisme de faisceaux inversibles sur $S$ dont le norme est 1 en toute
place, ce qui signifie $\|x_1/x_0\|^\sigma\geq 1$.

\begin{prop}\label{prop:rel-0}
Avec les notations pr\'ec\'edentes, si $\eps_Q:S\ra A^0$, il existe un point
$P\in\overline E(S)$ de hauteur relative nulle relevant $Q$ si et seulement si
le faisceau inversible m\'etris\'e $\eps_Q^*\mcL \in\hPic(S)$ est trivial.
\end{prop}
\begin{proof}
En effet, les consid\'erations qui pr\'ec\`edent montre que si $P$ est
un tel point, les deux sections $x_0:\mathscr O_S\ra \mathscr I$
et $x_1:\eps_Q^*\mcL^\vee\ra \mathscr I$ sont des isomorphismes de faisceaux
inversibles et $\|x_0/x_1\|^\sigma=1$ pour tout $\sigma\in
S(\C)$. 
Autrement dit, la section rationnelle $x_1\otimes x_0^{-1}$ de $\mcL$ est
sans p\^oles ni z\'eros et de norme~$1$. Cela signifie bien qu'elle
r\'ealise un isomorphisme $\mathscr O_S\simeq\mcL$ dans $\hPic(S)$.

Plus pr\'ecis\'ement, on voit que les points de hauteur
relative nulle de $\overline E(S)$ qui rel\`event la section $\eps_Q:S\ra A^0$
sont en bijection naturelle avec les trivialisations m\'etriques
du fibr\'e inversible m\'etris\'e $\eps_Q^*\mcL$.
\end{proof}

On d\'eduit de la proposition pr\'ec\'edente le fait 
suivant~\cite[proposition 2]{bertrand95}~: {\em si $K=\Q$ ou si $K$ est un corps
quadratique imaginaire, la nullit\'e de la hauteur de N\'eron--Tate
de $x\in A(K)$ relativement \`a $\mcL$ implique qu'il existe un point
de hauteur relative nulle dans l'extension param\'etr\'ee par
$\mcL$ relevant un multiple de $x$.} En effet, comme 
$\mcL|_{nx}=(\mcL|_x)^{\otimes n}$, on peut choisir $n$ de sorte
que $\mcL|_{nx}\simeq \mathscr O_S$ dans $\Pic(A)$~; 
il n'y a qu'une place archim\'edienne
et si $s$ est une base de $\mcL|_{nx}$, la nullit\'e de $\hdeg \mcL|_{nx}$
implique que $\|s\|=1$ et $\mcL|_{nx}=0$ dans $\hPic(S)$, ce qu'il fallait
d\'emontrer.

Consid\'erons alors un point $Q\in A^0(S)$ de la forme $[n]Q_1$,
avec $Q_1\in A^0(S)$ et $P\in E(S)$ un point de hauteur relative nulle 
qui rel\`eve $Q$, on peut consid\'erer un point $P_1$ relevant $Q_1$
tel que $[n]_E P_1=P$. Or, si le faisceau d'Arakelov $\eps_{Q_1}^*\mcL$
est de torsion, il n'est pas forc\'ement trivial. 
Gr\^ace au lemme suivant (dont il peut \^etre int\'eressant de remarquer que 
l'analogue g\'eom\'etrique est classique), 
cela signifie que $P_1$ n'est pas {\em a priori} d\'efini sur~$S$ :
\begin{lemm}
Soient $S$ le spectre de l'anneau des entiers d'un corps de nombres et
$(L,\|\cdot\|)$ un \'el\'ement de $\hPic(S)$, d'ordre fini.
Alors il existe une extension finie $f:S'\ra S$ telle que $f^*(L,\|\cdot\|)$
est nul dans $\hPic(S')$.
\end{lemm}
\begin{proof}
Soit $n\geq 1$ un entier tel que
$(L^{\otimes n},\|\cdot\|^n)=0$
et choisissons
une base $t_n\in L^{\otimes n}$ de norme 1 en toute place, soit
$\|t_n\|^\sigma=1$ pour tout $\sigma\in S(\C)$. Il existe une
extension finie $S'/S$ (par exemple,
l'anneau des entiers du corps de classe de Hilbert
du corps des fractions de $S$)
telle que $L'=L\otimes_{\mathscr O_S} \mathscr
O_{S'}$ est un $\mathscr O_{S'}$ module libre de rang~1~; soit donc
$s'$ une base de $L'$. Alors, $s^{\prime \otimes n}$ est une base
de $L^n \otimes \mathscr O_{S'}$ et il existe une unit\'e $u$ de $\mathscr
O_{S'}$ telle que $s^{\prime \otimes n}=u t_n$. Une extension $S''\ra S'$
telle que $u^{1/n}\in\mathscr O_{S''}$ permet de poser $s''=u^{-1/n} s'$~;
on constate que $s''$ est une base de $L''=L'\otimes_{\mathscr O_S}\mathscr
O_{S''}$ de norme 1 en toute place, si bien que $L\otimes \mathscr
O_{S''}=0$ dans $\hPic(S'')$.
\end{proof}

Montrons enfin comment les {\og points de Ribet\fg}
de~\cite{jacquinot-r87,bertrand96} s'interpr\`etent dans ce contexte, en 
supposant pour simplifier que
$A/S$ est un sch\'ema ab\'elien. 
\begin{prop}[{voir aussi~\protect\cite[th. 4]{bertrand96}}]
\label{prop:ribet-pt}
Supposons que $A$ est un sch\'ema ab\'elien sur~$S$.
Soient $f:A^\vee\ra A$ un $S$-morphisme de sch\'emas ab\'eliens et
$g=f-f^\vee: A^\vee_\eta\ra A_\eta$ qui est un endomorphisme antisym\'etrique.
Il existe au-dessus du point $g(\mcL)\in A(S)$ un $S$-point canonique
de hauteur relative nulle.
\end{prop}
\begin{proof}
Sur $X=A\times A^\vee$, consid\'erons la 
(bi)extension de Poincar\'e $\mathscr P_X$,
et de m\^eme sur $Y=A^\vee\times A$, m\'etris\'es de sorte
que le th\'eor\`eme du cube soit une isom\'etrie~\cite{moret-bailly85b}. 
Soit $s$ l'isomorphisme $Y\simeq X$ qui \'echange les facteurs~;
par unicit\'e du prolongement m\'etris\'e, 
on a un isomorphisme de faisceaux inversible m\'etris\'es $s^*\mathscr P_X
=\mathscr P_Y$ qui prolonge la bidualit\'e sur la fibre g\'en\'erique.

Le crit\`ere valuatif de propret\'e entra\^\i ne
que $g$ se prolonge en un unique endomorphisme $A^{\vee}\ra A$.
Alors, on a les \'egalit\'es entre faisceaux inversibles m\'etris\'es, 
qui r\'esultent de ce qu'elles sont vraies
sur $\eta$ et de l'unicit\'e du prolongement~: 
\begin{eqnarray*}
\mcL|_{f(\mcL)} &=& \mathscr P_X |_{f(\mcL),\mcL} 
= \mathscr P_Y |_{\mcL,f^\vee(\mcL)} \qquad\text{par dualit\'e}\\
&=& (s^*\mathscr P_X) |_{\mcL,f^\vee(\mcL)} 
=\mathscr P_X |_{f^\vee(\mcL),\mcL} 
= \mcL|_{f^\vee(\mcL)}, 
\end{eqnarray*}
si bien que $\mcL|_{g(\mcL)}$ est canoniquement
trivial, en tant que faisceau inversible 
m\'etris\'e sur~$S$. 
La preuve de la proposition~\ref{prop:rel-0}
montre que les $S$-points de hauteur 
relative nulle relevant $g(\mcL)$ correspondent aux isomorphismes
$\mcL|_{g(\mcL)}\ra (\mathscr O_S,\|1\|)$ dans $\hPic(S)$, d'o\`u un point
canonique d\'efini par l'isomorphisme ci-dessus.
\end{proof}
\begin{rema}
Prouvons que
le {\og point de Ribet\fg} d\'efini dans~\cite{jacquinot-r87} et consid\'er\'e 
dans~\cite{bertrand96} du point de vue des hauteurs est \'egal au point 
donn\'e par la proposition pr\'ec\'edente.
Consid\'erons comme dans~\cite[(4.1), p.~146]{jacquinot-r87} le diagramme
$$
\xymatrix{
  1 \ar[r]  & {\gm } \ar[r] \ar@{=}[d] & {f^*E} \ar[r] \ar[d]
       & {A^\vee} \ar[r] \ar_{f}[d]  & 0 \\
  1 \ar[r]  & {\gm } \ar[r] & {E} \ar[r] 
       & {A} \ar[r] & 0. }
$$
Choisissons $y\in f^*E (\eta)$ relevant $\mcL$, soit $x^1=f(y)\in E(\eta)$,
d'o\`u un 1--motif~$M_1\colon \Z \ra f^*E$.
Par d\'efinition de la dualit\'e de Cartier des 1--motifs,
le dual de~$M_1$ est un 1--motif~$M_2:\Z \ra E$. Soit $x^2\in E(\eta)$ l'image
de~$1$. Le point $x_f:=x^1-x^2\in E(\eta)$ rel\`eve $f(\mcL)-f^\vee(\mcL)$
et Bertrand prouve dans~\cite{bertrand96} que ce point est de hauteur relative 
nulle.  Quand $A/S$ est un sch\'ema ab\'elien, on peut faire cette 
construction sur~$S$ en choisissant $y\in f^*E(S)$ (si c'est possible)
et en raisonnant en termes de $S$--1--motifs.
Pour simplifier, raisonnons {\og localement sur $S\union\{\infty\}$\fg}
($S$ compactifi\'e en rajoutant les places \`a l'infini) --- le r\'esultat
\`a obtenir est de nature locale, en l'esp\`ece une hauteur locale
\`a calculer en un point ind\'ependant de~$y$ ---
et choisissons $y\in f^*E(S)$~; il correspond donc \`a
un isomorphisme $f^\vee(\mcL)|_\mcL \ra \mathscr O_S$ puisque $f^*E$ est
param\'etr\'ee par $f^\vee$.
Alors, la description sym\'etrique du 1--motif $M_1$ est la trivialisation
de $\mathscr P_Y|_{\mcL,f^\vee(\mcL)}=f^\vee(\mcL)|_{\mcL}$  que 
d\'efinit~$y$. Le 1--motif dual est alors donn\'e par la trivialisation
canonique de $\mathscr P_X|_{f(\mcL),\mcL}$ donn\'ee par la trivialisation
pr\'ec\'edente et la dualit\'e entre les $S$--sch\'emas ab\'eliens $A$ 
et $A^\vee$. Autrement dit, $x^2$ est donn\'e par la trivialisation
de $\mcL|_{f(\mcL)}$ qu'on en d\'eduit comme dans la preuve de la 
proposition et la diff\'erence $x^1-x^2$ est d\'efini par l'isomorphisme
$\mcL|_{f(\mcL)} \simeq \mcL|_{f^\vee(\mcL)}$ de la d\'emonstration
de la proposition~\ref{prop:ribet-pt}, isomorphisme qu'on a vu \^etre une 
isom\'etrie.
\end{rema}

\begin{rema} 
On v\'erifie ais\'ement que
les expressions explicites pour $\hdeg H_0$ et $\hdeg H_\infty$
que nous avons \'ecrites plus haut
donnent la d\'ecomposition de la hauteur relative en
une somme de hauteurs locales canoniques, comme dans~\cite{bertrand96,call-s93}.
\end{rema}

\section{M\'etriques ad\'eliques \emph{vs.} mod\`eles entiers}

On peut \'eviter les r\'ef\'erences aux mod\`eles de N\'eron
dans cet article, voire ne pas supposer que $A/S$ est un sch\'ema 
ab\'elien comme dans la proposition~\ref{prop:ribet-pt}
en faisant appel \`a la th\'eorie des m\'etriques ad\'eliques
due \`a S.~Zhang (cf.~\cite{zhang95b}). Soit $A_K$ une vari\'et\'e ab\'elienne
sur un corps de nombres~$K$ et $E_K$ une extension de $A_K$ par le
groupe multiplicatif $\gm$, donn\'ee par un faisceau inversible
$\mcL\in\Pic^ 0(A_K)$ et un {\og isomorphisme du carr\'e\fg}
$$ C(\mcL) := p_{12}^ *\mcL \otimes p_1^ *\mcL^ \vee \otimes p_2^
*\mcL^\vee \simeq \mathscr O_{A_K\times A_K}, $$
d'o\`u en particulier un isomorphisme canonique pour tout $n\in\Z$
$[n]^ *\mcL\simeq \mcL^ n$. Un entier $n\geq 2$ \'etant fix\'e,
le fibr\'e en droites $\mcL$ poss\`ede pour toute place~$v$ de $K$
une unique m\'etrique $v$-adique telle que cet isomorphisme
soit une isom\'etrie (th\'eor\`eme 2.2 de~\cite{zhang95b}
dont la d\'emonstration est une adaptation dans ce contexte
du proc\'ed\'e utilis\'e par Tate pour construire la hauteur normalis\'ee
sur une vari\'et\'e ab\'elienne).

De plus, cette m\'etrique $v$-adique rend le th\'eor\`eme du carr\'e une
isom\'etrie.
La m\'etrique $v$-adique sur $\mcL$ induit en effet une 
m\'etrique $v$-adique sur $C(\mathscr L)$ telle
que l'on ait une isom\'etrie
$[n]^*C(\mcL) \simeq C(\mcL^ n) \simeq C(\mcL)^ n$~; comme
$C(\mcL)$ est trivial, cette m\'etrique
$v$-adique est n\'ecessairement la m\'etrique triviale sur $\mathscr
O_{A_K\times A_K}$ telle que la section~$1$ a pour norme~$1$ en tout
point, cqfd.
En particulier, pour tout $p\in\Z$, l'isomorphisme
$[p]^*\mathscr L\simeq\mathscr L^p$ est une isom\'etrie.

Si $A_K$ a bonne r\'eduction en~$v$, 
l'unique mod\`ele entier $(A_v,\mcL_v)$,
$A_v$ \'etant un $\mathfrak o_v$ sch\'ema ab\'elien qui prolonge $A_K$
et $\mcL_v$ l'unique \'el\'ement de $\Pic^ 0(A_v) $ qui prolonge $A_K$,
nous fournit une m\'etrique $v$-adique canonique donnant en $x\in A_K(\bar
K_v)$ une norme $\leq 1$ si cette section est r\'eguli\`ere dans un voisinage
de l'adh\'erence de $x$ dans $A_v$. Cette m\'etrique co\"\i ncide
avec celle que nous avions d\'efinie auparavant.

La collection de ces m\'etriques $v$-adiques constitue de plus
une m\'etrique ad\'elique
au sens de {\em loc.\ cit.} Elles permettent par les m\^emes formules
que celles que nous avons donn\'ees dans le cas archim\'edien
de construire une m\'etrique ad\'elique canonique sur le faisceau
$\mathscr O_\P(1)$ de la compactification
$\P(\mathscr O_{A_K}\oplus\mcL^ \vee)$,
et cette m\'etrique ad\'elique v\'erifie
des propri\'et\'es tout \`a fait analogues
\`a celles que nous avons \'etablies dans cet article.

Explicitons maintenant ce qu'est une m\'etrique 
ad\'elique pour un fibr\'e en droites sur $\Spec K$.
Un tel fibr\'e est un $K$-espace vectoriel $L$ de
dimension~$1$ ; fixons une base $e$ de $L$.  Une m\'etrique ad\'elique
sur $L$ revient alors \`a une collection $(\norm{e}_v)\in\prod_v |K^*|_v$
pour toutes les places $v$ de $K$, $|K^*|_v$ d\'esignant le groupe
des normes $v$-adiques des \'el\'ements non nuls de $K$, telle que
$\norm{e}_v\leq 1$ pour presque tout~$v$.
Il poss\`ede un degr\'e arithm\'etique, d\'efini par
$$\hdeg (L,(\norm{\cdot}_v)) = -\sum_v \log_v \norm{e}_v, $$
o\`u $\log_v:|K|_v^*\ra\Q$ est la normalisation naturelle du logarithme
de sorte que la formule du produit s'écrit $\sum_v \log_v |x|_v=0$ pour tout
$x\in K^ *$. Ainsi, le degr\'e arithm\'etique ne d\'epend pas du choix de la
base. En fait, cette description rend apparent
qu'il existe un unique isomorphisme
du groupe des classes d'isomorphismes de fibr\'es en droites sur $\Spec K$
avec m\'etriques ad\'eliques sur le groupe $\hPic (\Spec\mathfrak o_K)$
compatible aux degr\'es arithm\'etiques.

Cette notion de m\'etriques ad\'eliques
fournit ainsi une construction alternative des hauteurs
canoniques sur l'extension $E_K$.
L'analogue de la proposition~\ref{prop:rel-0} est que les points de $E_K(K)$
de hauteur relative nulle relevant un point $x\in A_K(K)$ sont
en bijection naturelle avec les trivialisations isom\'etriques
de $\mcL_x$. La proposition~\ref{prop:ribet-pt} se prouve dans ce cadre
en rempla\c{c}ant les isomorphismes dans les groupes de Picard
compactifi\'es  utilis\'es dans la preuve
de cette proposition par des isom\'etries de fibr\'es en droites munis de
m\'etriques ad\'eliques. Ainsi formulée, elle s'étend au cas de mauvaise
réduction.

\def\No{\leavevmode\kern-.1em\raise.7ex\hbox{\smaller[4]\upshape
  o}\kern.17em\relax}%

\end{document}